\documentclass[11pt]{article}
\usepackage{epic,latexsym,amssymb}

\newtheorem{thm}{Theorem}
\newtheorem{lem}[thm]{Lemma}

\newtheorem{ob}{Observation}
\newtheorem{prop}{Proposition}

\newcommand{\diam}{{\rm diam}}

\newcommand{\ad}{{\rm ad}}
\newcommand{\mad}{{\rm mad}}

\newcommand{\cP}{{\cal P}}

\newcommand{\barG}{\overline{G}}
\newcommand{\barK}{\overline{K}}

\newcommand{\TR}[1]{\mbox{$\tau(#1)$}}

\newcommand{\TRr}[1]{\mbox{$\tau_r(#1)$}}
\newcommand{\TRo}[1]{\mbox{$\tau_1(#1)$}}

\newcommand{\1}{ \vspace{0.1cm} }

\newcommand{\mop}{{\rm mop}}
\newcommand{\Mop}{{\rm Mop}}
\newcommand{\mpp}{{\rm mp}}
\newcommand{\Mpp}{{\rm Mp}}

\newcommand{\cG}{{\cal G}}
\newcommand{\cR}{{\cal R}}

\newcommand{\cPP}{{\cal P}}
\newcommand{\cOP}{{\cal OP}}

\newenvironment{unnumbered}[1]{\trivlist \item [\hskip \labelsep {\bf
#1}]\ignorespaces\it}{\endtrivlist}

\def\vertex(#1){\put(#1){\circle*{2}}}
\def\vertexo(#1){\put(#1){\circle{2}}}
\def\vert(#1){\put(#1){\circle*{1.5}}}
\def\verto(#1){\put(#1){\circle{1.5}}}
\def\lab(#1)#2{\put(#1){\makebox(0,0)[c]{#2}}}
\begin{document}

\title{Directed Domination in Oriented Graphs}

\author{$^1$Yair Caro and $^2$Michael A. Henning\thanks{Research supported in part by the South African National Research Foundation} \\
\\
$^1$Department of Mathematics and Physics\\
University of Haifa-Oranim \\
Tivon 36006, Israel \\
Email: yacaro@kvgeva.org.il \\
\\
$^2$Department of Mathematics \\
University of Johannesburg \\
Auckland Park 2006, South Africa \\
Email: mahenning@uj.ac.za }

\date{}
\maketitle

\begin{abstract}
A directed dominating set in a directed graph $D$ is a set $S$ of
vertices of $V$ such that every vertex $u \in V(D) \setminus S$ has
an adjacent vertex $v$ in $S$ with $v$ directed to $u$. The directed
domination number of $D$, denoted by $\gamma(D)$, is the minimum
cardinality of a directed dominating set in $D$. The directed
domination number of a graph $G$, denoted $\Gamma_d(G)$, which is the
maximum directed domination number $\gamma(D)$ over all orientations
$D$ of $G$. The directed domination number of a complete graph was
first studied by Erd\"{o}s [Math. Gaz. 47 (1963), 220--222], albeit
in disguised form. We extend this notion to directed domination of
all graphs. If $\alpha$ denotes the independence number of a graph
$G$, we show that if $G$ is a bipartite graph, we show that
$\Gamma_d(G) = \alpha$. We present several lower and upper bounds on
the directed domination number. \end{abstract}


\newpage
\section{Introduction}

An \emph{asymmetric digraph} or \emph{oriented graph} $D$ is a
digraph that can be obtained from a graph $G$ by assigning a
direction to (that is, orienting) each edge of $G$. The resulting
digraph $D$ is called an \emph{orientation} of $G$. Thus if $D$ is
an oriented graph, then for every pair $u$ and $v$ of distinct
vertices of $D$, at most one of $(u,v)$ and $(v,u)$ is an arc of
$D$.
A \emph{directed dominating set}, abbreviated DDS, in a directed
graph $D = (V,A)$ is a set $S$ of vertices of $V$ such that every
vertex in $V \setminus S$ is dominated by some vertex of $S$; that
is, every vertex $u \in V \setminus S$ has an adjacent vertex $v$ in
$S$ with $v$ directed to $u$. Every digraph has a DDS since the
entire vertex set of the digraph is such a set. The \emph{directed
domination number} of a \emph{directed graph} $D$, denoted by
$\gamma(D)$, is the minimum cardinality of a DDS in $D$. A DDS of
$D$ of cardinality~$\gamma(D)$ is called a $\gamma(D)$-set. Directed
domination in digraphs is well studied
(cf.~\cite{ArJaVo07,BhVi05,cdss,ChHaYu97,ChVaYu96,Fu68,GhLaPi98,hhs2,Le98,ReMcHeHe04}).

We define the \emph{lower directed domination number} of a
\emph{graph} $G$, denote $\gamma_d(G)$, to be the minimum directed
domination number $\gamma(D)$ over all orientations $D$ of $G$; that
is,
\[
\gamma_d(G) =  \min \{ \gamma(D) \mid \mbox{ over all orientations
$D$ of $G$} \}.
\]
The \emph{upper directed domination number}, or simply the
\emph{directed domination number}, of a \emph{graph} $G$, denoted
$\Gamma_d(G)$, is defined as the maximum directed domination number
$\gamma(D)$ over all orientations $D$ of $G$; that is,
\[
\Gamma_d(G) = \max \{ \gamma(D) \mid \mbox{ over all orientations
$D$ of $G$} \}.
\]

\subsection{Motivation}

The directed domination number of a complete graph was first studied
by Erd\"{o}s~\cite{Er63} albeit in disguised form. In 1962,
Sch\"{u}tte~\cite{Er63} raised the question of given any positive
integer $k > 0$, does there exist a tournament $T_{n(k)}$ on $n(k)$
vertices in which for any set $S$ of $k$ vertices, there is a vertex
$u$ which dominates all vertices in $S$. Erd\"{o}s~\cite{Er63}
showed, by probabilistic arguments, that such a tournament
$T_{n(k)}$ does exist, for every positive integer $k$. The proof of
the following bounds on the directed domination number of a complete
graph are along identical lines to that presented by
Erd\"{o}s~\cite{Er63}. This result can also be found
in~\cite{ReMcHeHe04}. Throughout this paper, $\log$ is to the
base~$2$ while $\ln$ denotes the logarithm in the natural base $e$.

\begin{thm}{\rm (Erd\"{o}s~\cite{Er63})}
For every integer $n \ge 2$, $\log n - 2\log (\log n) \le
\Gamma_d(K_n) \le \log (n+1)$. \label{t:KnLow}
\end{thm}

In this paper, we extend this notion of directed domination in a
complete graph to directed domination of all graphs.

\subsection{Notation}

For notation and graph theory terminology we in general
follow~\cite{hhs1}. Specifically, let $G = (V, E)$ be a graph with
vertex set $V$ of order~$n = |V|$ and edge set $E$ of size~$m =
|E|$, and let $v$ be a vertex in $V$. The \emph{open neighborhood}
of $v$ is $N_G(v) = \{u \in V \, | \, uv \in E\}$ and the
\emph{closed neighborhood of $v$} is $N_G[v] = \{v\} \cup N_G(v)$.
If the graph $G$ is clear from context, we simply write $N(v)$ and
$N[v]$ rather than $N_G(v)$ and $N_G[v]$, respectively. For a set $S
\subseteq V$, the subgraph induced by $S$ is denoted by $G[S]$. If
$A$ and $B$ are subsets of $V(G)$, we let $[A,B]$ denote the set of
all edges between $A$ and $B$ in $G$. We denote the diameter of $G$
by $\diam(G)$.

We denote the \emph{degree} of $v$ in $G$ by $d_G(v)$, or simply by
$d(v)$ if the graph $G$ is clear from context. The minimum degree
among the vertices of $G$ is denoted by $\delta(G)$, and the maximum
degree by~$\Delta(G)$.
The \emph{maximum average degree} in $G$, denoted by $\mad(G)$, is
defined as the maximum of the average degrees $\ad(H) = 2
|E(H)|/|V(H)|$ taken over all subgraphs $H$ of $G$.

The parameter $\gamma(G)$ denotes the \emph{domination number} of
$G$. The parameters $\alpha(G)$ and $\alpha'(G)$ denote the (vertex)
\emph{independence number} and the \emph{matching number},
respectively, of $G$, while $\chi(G)$ and $\chi'(G)$ denote the
\emph{chromatic number} and \emph{edge chromatic number},
respectively, of $G$. The \emph{covering number} of $G$, denoted by
$\beta(G)$,  is the minimum number vertices that covers all the edges
of $G$. The \emph{clique number} of $G$, denoted by $\omega(G)$, is
the maximum cardinality of a clique in $G$.

A vertex $v$ in a digraph $D$  \emph{out-dominates}, or simply
\emph{dominates}, itself as well as all vertices $u$ such that $(v,
u)$ is an arc of $D$. The \emph{out-neighborhood} of $v$, denoted
$N^+(v)$, is the set of all vertices $u$ adjacent from $v$ in $D$;
that is, $N^+(v) = \{u \mid (v, u) \in A(D)\}$. The
\emph{out-degree} of $v$ is given by $d^+(v) = |N^+(v)|$, and the
maximum out-degree among the vertices of $D$ is denoted by
$\Delta^+(D)$. The \emph{in-neighborhood} of $v$, denoted $N^-(v)$,
is the set of all vertices $u$ adjacent to $v$ in $D$; that is,
$N^-(v) = \{u \mid (u,v) \in A(D)\}$. The \emph{in-degree} of $v$ is
given by $d^-(v) = |N^-(v)|$. The \emph{closed in-neighborhood} of
$v$ is the set $N^-[v] = N^-(v) \cup \{v\}$. The maximum in-degree
among the vertices of $D$ is denoted by $\Delta^-(D)$.

A \emph{hypergraph} $H = (V,E)$ is a finite set $V$ of elements,
called \emph{vertices}, together with a finite multiset $E$ of
subsets of $V$, called \emph{edges}. A $k$-\emph{edge} in $H$ is an
edge of size~$k$.
The hypergraph $H$ is said to be $k$-\emph{uniform} if every edge of
$H$ is a $k$-edge. A subset $T$ of vertices in a hypergraph $H$ is a
\emph{transversal} (also called \emph{vertex cover} or \emph{hitting
set} in many papers) if $T$ has a nonempty intersection with every
edge of $H$. The \emph{transversal number} $\TR{H}$ of $H$ is the
minimum size of a transversal in $H$.
For a digraph $D = (V,E)$, we denote by $H_D$ the \emph{closed
in-neighborhood hypergraph}, abbreviated CINH, of $D$; that is, $H_D
= (V,C)$ is the hypergraph with vertex set $V$ and with edge set $C$
consisting of the closed in-neighborhoods of vertices of $V$ in $D$.

\section{Observations}

We show first that the lower directed domination number of a graph is
precisely its domination number.

\begin{ob}
For every graph $G$, $\gamma_d(G) = \gamma(G)$. \label{o:lower}
\end{ob}
\textbf{Proof.} Let $S$ be a $\gamma(G)$-set and let $D$ be an
orientation obtained from $G$ by directing all edges in $[S,V
\setminus S]$ from $S$ to $V \setminus S$ and directing all other
edges arbitrarily. Then, $S$ is a DDS of $D$, and so $\gamma_d(G)
\le \gamma(D) \le |S| = \gamma(G)$. However if $D$ is an orientation
of a graph $G$ such that $\gamma_d(G) = \gamma(D)$, and if $S$ is a
$\gamma(D)$-set, then $S$ is also a dominating set of $G$, and so
$\gamma(G) \le |S| = \gamma_d(G)$. Consequently, $\gamma_d(G) =
\gamma(G)$.~$\Box$

\medskip
In view of Observation~\ref{o:lower}, it is not interesting to ask
about the lower directed domination number, $\gamma_d(G)$, of a graph
$G$ since this is precisely its domination number, $\gamma(G)$, which
is very well studied. We therefore focus our attention on the (upper)
directed domination number of a graph. As a consequence of
Theorem~\ref{t:KnLow}, we establish a lower bound on the directed
domination number of an arbitrary graph.

\begin{ob}
For every graph $G$ on $n$ vertices, $\Gamma_d(G) \ge \log n - 2\log
(\log n)$. \label{oblow}
\end{ob}
\textbf{Proof.}  Let $D$ be an orientation of the edges of a complete
graph $K_n$ on the same vertex set as $G$ such that $\Gamma_d(K_n) =
\gamma(D)$. Let $D_G$ be the orientation of $D$ induced by arcs of
$D$ corresponding to edges of $G$. Then, $\Gamma_d(G) \ge \gamma(D_G)
\ge \gamma(D) = \Gamma_d(K_n)$. The desired lower bound now follows
from Theorem~\ref{t:KnLow}.~$\Box$

\begin{ob}
If $H$ is an induced subgraph of a graph $G$, then $\Gamma_d(G) \ge
\Gamma_d(H)$. \label{o:induced}
\end{ob}
\textbf{Proof.} Let $G = (V,E)$ and let $U = V(H)$. Let $D_H$ be an
orientation of $H$ such that $\Gamma_d(H) = \gamma(D_H)$. We now
extend the orientation $D_H$ of $H$ to an orientation $D$ of $G$ by
directing all edges in $[U,V \setminus U]$ from $U$ to $V \setminus
U$ and directing all edges with both ends in $V \setminus U$
arbitrarily. Then, $\Gamma_d(G) \ge \gamma(D) \ge \gamma(D_H) =
\Gamma_d(H)$.~$\Box$

\begin{ob}
If $H$ is a spanning subgraph of a graph $G$, then $\Gamma_d(G) \le
\Gamma_d(H)$. \label{o:spanning}
\end{ob}
\textbf{Proof.} Let $D$ be an arbitrary orientation of $G$, and let
$D_H$ be the orientation of $H$ induced by $D$. Since adding arcs
cannot increase the directed domination number, we have that
$\gamma(D) \le \gamma(D_H)$. This is true for every orientation of
$G$. Hence, $\Gamma_d(G) \le \Gamma_d(H)$.~$\Box$

\medskip
Hakimi~\cite{Ha65} proved that a graph $G$ has an orientation $D$
such that $\Delta^+(D) \le k$ if and only if $\mad(G) \le 2k$. This
implies the following result.

\begin{ob}{\rm (\cite{Ha65})}
Every graph $G$ has an orientation $D$ such that $\Delta^+(D) \le
\lceil \mad(G)/2 \rceil$.
 \label{o:Hakimi}
\end{ob}

\section{Bounds}

In this section, we establish bounds on the directed domination
number of a graph. We first present lower bounds on the directed
domination number of a graph.

\begin{thm}
Let $G$ be a graph of order~$n$. Then the following holds. \\
{\rm (a)} $\Gamma_d(G) \ge \alpha(G) \ge \gamma(G)$. \\
{\rm (b)} $\Gamma_d(G) \ge n/\chi(G)$. \\
{\rm (c)} $\Gamma_d(G) \ge \lceil (\diam(G)+1)/2) \rceil$. \\
{\rm (d)}
$\Gamma_d(G) \ge n/( \lceil \mad(G)/2 \rceil + 1)$. \\
\label{t:lowerbd}
\end{thm}
\textbf{Proof.} Since every maximal independent set in a graph is a
dominating set in the graph, we recall that $\gamma(G) \le \alpha(G)$
holds for every graph $G$. To prove that $\alpha(G) \le \Gamma_d(G)$,
let $A$ be a maximum independent set in $G$ and let $D$ be the
digraph obtained from $G$ by orienting all arcs from $A$ to $V
\setminus A$ and orienting all arcs in $G[V \setminus A]$, if any,
arbitrarily. Since every DDS of $D$ contains $A$, we have $\gamma(D)
\ge |A|$. However the set $A$ itself is a DDS of $D$, and so
$\gamma(D) \le |A|$. Consequently, $\Gamma_d(G) \ge \gamma(D) = |A| =
\alpha(G)$. This establishes Part~(a).
Parts~(b) and~(c) follows readily from Part~(a) and the observations
that $\alpha(G) \ge n/\chi(G)$ and $\alpha(G) \ge \lceil
(\diam(G)+1)/2) \rceil$.
%
By Observations~\ref{o:Hakimi}, there is an orientation $D$ of $G$
such that $\Delta^+(D) \le \lceil \mad(G)/2 \rceil$. Let $S$ be a
$\gamma(D)$-set. Then, $\displaystyle{ V \setminus S \subseteq
\cup_{v \in S} N^+(v) }$, and so $n - |S| = |V \setminus S| \le
\sum_{v \in S} d^+(v) \le |S| \cdot \Delta^+(D)$, whence $\gamma(D) =
|S| \ge n/(\Delta^+(D) + 1) \ge n/( \lceil \mad(G)/2 \rceil + 1)$.
This establishes Part~(d).~$\Box$

\medskip
We remark that since $\mad(G) \le \Delta(G)$ for every graph $G$, as
an immediate consequence of Theorem~\ref{t:lowerbd}(d) we have that
$\Gamma_d(G) \ge n/( \lceil \Delta(G)/2 \rceil + 1)$.

\medskip
Next we consider upper bounds on the directed domination number of a
graph. The following lemma will prove to be useful.

\begin{lem}
Let $G = (V,E)$ be a graph and let $V_1,V_2,\ldots,V_k$ be subsets
of $V$, not necessarily disjoint, such that $\cup_{i=1}^k V_i =
V(G)$. For $i = 1,2,\ldots,k$, let $G_i = G[V_i]$. Then,
\[
\Gamma_d(G) \le \sum_{i=1}^k \Gamma_d(G_i).
\]
 \label{usefulLem}
\end{lem}
\textbf{Proof.} Consider an arbitrary orientation $D$ of $G$. For
each $i =1,2,\ldots,k$, let $D_i$ be the orientation of the edges of
$G_i$ induced by $D$ and let $S_i$ be a $\gamma(D_i)$-set. Then,
$\Gamma_d(G_i) \ge \gamma(D_i) = |S_i|$ for each $i$. Since the set
$S = \cup_{i=1}^k S_i$ is a DDS of $D$, we have that $\gamma(D) \le
|S| \le \sum_{i=1}^k |S_i| \le \sum_{i=1}^k \Gamma_d(G_i)$.
%
Since this is true for every orientation $D$ of $G$, the desired
upper bound on $\Gamma_d(G)$ follows.~$\Box$

\medskip
As a consequence of Lemma~\ref{usefulLem}, we have the following
upper bounds on the directed domination number of a graph.

\begin{thm}
Let $G$ be a graph of order~$n$. Then the following holds. \\
{\rm (a)} $\Gamma_d(G) \le n - \alpha'(G)$. \\
{\rm (b)} If $G$ has a perfect matching, then $\Gamma_d(G) \le n/2$.
\\
{\rm (c)} $\Gamma_d(G) \le n$ with equality if and only if $G =
\overline{K}_n$. \\
{\rm (d)} If $G$ has minimum degree~$\delta$ and $n \ge 2\delta$,
then $\Gamma_d(G) \le n - \delta$. \\
{\rm (e)} $\Gamma_d(G) = n -1$ if and only if every component of $G$
is a $K_1$-component, except for one \\ \hspace*{0.5cm} component
which is either a star or a complete graph $K_3$.
 \label{t:uppbd}
\end{thm}
\textbf{Proof.} (a) Let $M = \{u_1v_1,u_2v_2, \ldots, u_tv_t\}$ be a
maximum matching in $G$, and so $t = \alpha'(G)$. For $i =
1,2,\ldots,t$, let $V_i = \{u_i,v_i\}$. If $n > 2t$, let
$(V_{t+1},\ldots,V_{n-2t})$ be a partition of the remaining vertices
of $G$ into $n-2t$ subsets each consisting of a single vertex. By
Lemma~\ref{usefulLem}, $\Gamma_d(G) \le \sum_{i=1}^n \Gamma_d(G_i) =
t + (n-2t) = n - t = n - \alpha'(G)$. Part~(b) is an immediate
consequence of Part~(a). Part~(c) is an immediate consequence of
Part~(a) and the observation that $\alpha'(G) = 0$ if and only if $G
= \overline{K}_n$.

(d) It is well known (see, for example, Bollob\'{a}s~\cite{Bo04},
pp. 87) that if $G$ has $n$ vertices and minimum degree~$\delta$
with $n \ge 2\delta$, then $\alpha'(G) \ge \delta$. Hence by
Part~(a) above, $\Gamma_d(G) \le n - \delta$.

(e) Suppose that $\Gamma_d(G) = n - 1$. Then by Part~(a) above,
$\alpha'(G) = 1$. However every connected graph $F$ with $\alpha'(F)
= 1$ is either a star or a complete graph $K_3$. Hence, either $G$
is the vertex disjoint union of a star and isolated vertices or of a
complete graph $K_3$ and isolated vertices.~$\Box$

\medskip
We establish next that the directed domination number of a bipartite
graph is precisely its independence number. For this purpose, recall
that K\"{o}nig~\cite{Ko31} and Egerv\'{a}ry~\cite{Eg31} showed that
if $G$ is a bipartite graph, then $\alpha'(G) = \beta(G)$. Hence by
Gallai's Theorem~\cite{Ga59}, if $G$ is a bipartite graph of
order~$n$, then $\alpha(G) + \alpha'(G) = n$.

\begin{thm}
If $G$ is a bipartite graph, then $\Gamma_d(G) = \alpha(G)$.
\label{bipart}
\end{thm}
\textbf{Proof.} Since $G$ is a bipartite graph, we have that $n -
\alpha'(G) = \alpha(G)$. Thus by Theorem~\ref{t:lowerbd}(a) and
Theorem~\ref{t:uppbd}(b), we have that $\alpha(G) \le \Gamma_d(G)
\le n - \alpha'(G) = \alpha(G)$. Consequently, we must have equality
throughout this inequality chain. In particular, $\Gamma_d(G) =
\alpha(G)$.~$\Box$

\section{Relation to other Parameters}

The following result establishes an upper bound on the directed
domination of a graph in terms of its independence number and
chromatic number.

\newpage
\begin{thm}
For every graph $G$, we have $\Gamma_d(G) \le \alpha(G) \cdot \lceil
\chi(G)/2 \rceil$. \label{indepcolor}
\end{thm}
\textbf{Proof.} Let $G$ have order~$n$. If $\chi(G) = 1$, then $G$
is the empty graph, $\overline{K}_n$ and so $\Gamma_d(G) = n =
\alpha(G)$, while if $\chi(G) = 2$, then $G$ is a bipartite graph,
and so by Theorem~\ref{bipart}, $\Gamma_d(G) = \alpha(G)$. In both
cases, $\alpha(G) = \alpha(G) \cdot \lceil \chi(G)/2 \rceil$, and so
$\Gamma_d(G) = \alpha(G) \cdot \lceil \chi(G)/2 \rceil$. Hence we
may assume that $\chi(G) \ge 3$.
If $\chi(G) = 2k$ for some integer $k \ge 2$, then let
$V_1,V_2,\ldots,V_{2k}$ denote the color classes of $G$. For $i =
1,2,\ldots,k$, let $G_i$ be the subgraph $G[V_{2i-1} \cup V_{2i}]$
of $G$ induced by $V_{2i-1}$ and $V_{2i}$ and note that $G_i$ is a
bipartite graph. By Theorem~\ref{bipart}, $\Gamma_d(G_i) =
\alpha(G_i) \le \alpha(G)$ for all $1,2,\ldots,k$. Hence by
Lemma~\ref{usefulLem}, $\Gamma_d(G) \le \sum_{i=1}^k \Gamma_d(G_i)
\le k \alpha(G) = \alpha(G) \cdot \lceil \chi(G)/2 \rceil$, as
desired.
If $\chi(G) = 2k+1$ for some integer $k \ge 1$, then let
$V_1,V_2,\ldots,V_{2k+1}$ denote the color classes of $G$. For $i =
1,2,\ldots,k$, let $H_i$ be the subgraph of $G$ induced by
$V_{2i-1}$ and $V_{2i}$ and note that $H_i$ is a bipartite graph.
Further let $H_{k+1} = G[V_{2k+1}]$, and so $H_{k+1}$ is an empty
graph on $|V_{2k+1}| \le \alpha(G)$ vertices. By
Lemma~\ref{usefulLem}, $\Gamma_d(G) \le \sum_{i=1}^{k+1}
\Gamma_d(H_i) \le (k+1) \alpha(G) = \alpha(G) \cdot \lceil \chi(G)/2
\rceil$.~$\Box$

\medskip
As shown in the proof of Theorem~\ref{indepcolor}, the upper bound
of Theorem~\ref{indepcolor} is always attained if $\chi(G) \le 2$.
We remark that if $\chi(G) = 3$ or $\chi(G) = 4$, then the upper
bound of Theorem~\ref{indepcolor} is achievable by taking, for
example, $G = rK_t$ where $t \in \{3,4\}$ and $r$ is some positive
integer. In this case, $\chi(G) = t$ and $\Gamma_d(G) = 2r =
\alpha(G) \cdot \lceil \chi(G)/2 \rceil$.

\begin{thm}
If $G$ is a graph of order~$n$, then $\Gamma_d(G) \le n - \lfloor
\chi(G)/2 \rfloor$. \label{colorB}
\end{thm}
\textbf{Proof.} If $\chi(G) = 1$, then the bound is immediate since
$\Gamma_d(G) \le n$ by Theorem~\ref{t:uppbd}(c). Hence we may assume
that $\chi(G) = k \ge 2$. Let $V_1,V_2,\ldots,V_{k}$ denote the
color classes of $G$. By the minimality of the coloring, there is an
edge between every two color classes. In particular for $i =
1,2,\ldots,\lfloor k/2 \rfloor$, there is an edge between $V_{2i-1}$
and $V_{2i}$, and so $\alpha'(G) \ge \lfloor k/2 \rfloor$. Hence by
Theorem~\ref{t:uppbd}(a), $\Gamma_d(G) \le n - \alpha'(G) \le n -
\lfloor k/2 \rfloor$.~$\Box$

\medskip
We remark that the bound of Theorem~\ref{colorB} is achievable for
graphs with small chromatic number as may be seen by considering the
graph $G = \barK_{n-k} \cup K_k$ where $1 \le k \le 4$ and $n > k$.
We show next that the directed domination of a graph is at most the
average of its order and independence number. For this purpose, we
recall the Gallai-Milgram Theorem~\cite{GaMi60} for oriented graphs
which states that in every oriented graph $G = (V,E)$, there is a
partition of $V$ into at most $\alpha(G)$ vertex disjoint directed
paths.

\begin{thm}
If $G$ is a graph of order~$n$, then $\Gamma_d(G) \le (n +
\alpha(G))/2$. \label{alpha}
\end{thm}
\textbf{Proof.} Let $D$ be an orientation of $G$. By the
Gallai-Milgram Theorem for oriented graphs, there is a partition
$\cP  = \{P_1,P_2, \ldots, P_t\}$ of $V(D)$ into $t$ vertex disjoint
directed paths where $t \le \alpha(G)$. For $i = 1,2,\ldots,t$, let
$|P_i| = p_i$, and so $\sum_{i=1}^t p_i = n$. By
Lemma~\ref{usefulLem}, $\Gamma_d(G) \le \sum_{i=1}^t \Gamma_d(P_i) =
\sum_{i=1}^t \lceil p_i/2 \rceil \le \sum_{i=1}^t (p_i+1)/2 =
(\sum_{i=1}^t p_i/2) + t/2 = (n + \alpha(G))/2$.~$\Box$

\medskip
That the bound of Theorem~\ref{alpha} is best possible, may be seen
by considering, for example, the graph $G = rK_3 \cup sK_1$ of
order~$n = 3r + s$ with $\alpha(G) = r + s$ and $\Gamma_d(G) = 2r +
s = (n + \alpha(G))/2$.

The following result establishes an upper bound on the directed
domination of a graph in terms of the chromatic number of its
complement.

\begin{thm}
If $G$ is a graph of order~$n$, then $\displaystyle{ \Gamma_d(G) \le
\chi(\barG) \cdot \log \left( \left\lceil \frac{n}{\chi(\barG)}
\right\rceil + 1 \right) }$.
\label{color}
\end{thm}
\textbf{Proof.} Let $t = \chi(\barG)$ and consider a
$\chi(\barG)$-coloring of the complement $\barG$ of $G$ into $t$
color classes $Q_1, Q_2, \ldots, Q_t$, where $|Q_i| = q_i$ for $i =
1,2, \ldots, t$. For each $i = 1,2, \ldots, t$, the subgraph
$G[Q_i]$ of $G$ induced by $Q_i$ is a clique. We now consider an
arbitrary orientation $D$ of $G$, and we let $D_i = D[Q_i]$ denote
the orientation of the edges of the clique $G[Q_i]$ induced by $D$.
Then,
\[
\gamma(D) \le \sum_{i=1}^t \gamma(D_i) \le \sum_{i=1}^t
\Gamma_d(Q_i) = \sum_{i=1}^t \Gamma_d(K_{q_i}).
\]
This is true for every orientation $D$ of $G$, and so, by
Theorem~\ref{t:KnLow}, we have that $\Gamma_d(G) \le \sum_{i=1}^t
\log (q_i + 1)$,
where $\sum_{i=1}^t q_i = n$. By convexity the right hand side
attains its maximum when all summands are as equal as possible; that
is, some of the summands are $\lfloor n/t \rfloor$ and some are
$\lceil n/t \rceil$. Hence, $\Gamma_d(G) \le t \log (\lceil n/t
\rceil + 1)$.~$\Box$

\medskip
As a consequence of Theorem~\ref{color}, we have the following
result on the directed domination number of a dense graph with large
minimum degree.

\begin{thm}
If $G$ is a graph on $n$ vertices with minimum degree $\delta(G) \ge
(k-1)n/k$ where $k$ divides $n$, then $\Gamma_d(G) \le n \log
(k+1)/k$. \label{dense}
\end{thm}
\textbf{Proof.} Since $k \, | \, n$, we note that $n = k t$ and
$\delta(G) \ge (k - 1) t$ for some integer~$t$. By the well-known
Hajnal-Szemer\'{e}di Theorem~\cite{HaSz70}, the graph $G$ contains
$t$ vertex disjoint copies of $K_k$. Further, $\chi(\barG) \le t$.
Thus applying Theorem~\ref{color}, we have that $\Gamma_d(G) \le t
\log (k + 1)  = n \log (k + 1) / k$.~$\Box$

\section{Special Families of Graphs}

In this section, we consider the (upper) directed domination number
of special families of graph. As remarked earlier, the directed
domination number of a complete graph $K_n$ is determined by
Erd\"{o}s~\cite{Er63} in Theorem~\ref{t:KnLow}, while the directed
domination number of a bipartite graph is precisely its independence
number (see Theorem~\ref{bipart}).

\subsection{Regular Graphs}

For each given $\delta \ge 1$, applying Theorem~\ref{t:lowerbd}(a)
to the graph $G = K_{\delta,n - \delta}$ yields $\Gamma_d(G) \ge n -
\delta$. Hence without regularity, we observe that for each fixed
$\delta \ge 1$, there exists a graph $G$ of order~$n$ and minimum
degree~$\delta$ satisfying $\Gamma_d(G) \ge n - \delta$.
With regularity, the directed domination number of a graph may be
much smaller. For a given $r$, let $n = k(r+1)$ for some integer~$k$
and let $G$ consist of the disjoint union of $k$ copies of
$K_{r+1}$. Let $G_1, G_2, \ldots, G_k$ denote the components of $G$.
Each component of $G$ is $r$-regular, and by Theorem~\ref{t:KnLow},
$\Gamma_d(G) = \sum_{i=1}^k \Gamma_d(G_i) = \sum_{i=1}^k
\Gamma_d(K_{r+1}) \le k \log (r+2) = n \log (r+2) / (r+1). $ Hence
there exist $r$-regular graphs of order~$n$ with $ \Gamma_d(G) \le n
\log (r+2) / (r+1)$. In view of these observations it is of interest
to investigate the directed domination number of regular graphs.

In 1964, Vizing proved his important edge-coloring result which
states that every graph $G$ satisfies $\Delta(G) \le \chi'(G) \le
\Delta(G) + 1$. As a consequence of Vizing's Theorem, we have the
following upper bound on the directed domination number of a regular
graph.

\begin{thm}
For $r \ge 2$, if $G$ is an $r$-regular graph of order~$n$, then
\[
\Gamma_d(G) \le n(r+2)/2(r+1).
\]
\label{t:regular1}
\end{thm}
\textbf{Proof.} By Vizing's Theorem, $\chi'(G) \le r + 1$. Consider
an edge coloring of $G$ using $\chi'(G)$-colors. The edges in each
color class form a matching in $G$, and so the matching number of
$G$ is at least the size of a largest color class in $G$. Hence if
$G$ has size~$m$, we have $\alpha'(G) \ge m/\chi'(G) \ge m/(r+1) =
nr/2(r+1)$. Hence by Theorem~\ref{t:uppbd}(a), $\Gamma_d(G) \le n -
\alpha'(G) \le n - nr/2(r+1) = n(r+2)/2(r+1)$.~$\Box$

\medskip
As a special case of Theorem~\ref{t:regular1}, we have that
$\Gamma_d(G) \le 2n/3$ if $G$ is a $2$-regular graph. We next
characterize when equality is achieved in this bound.

\begin{prop}
Let $G$ be a $2$-regular graph on $n \ge 3$ vertices. Then the
following holds. \\
\indent {\rm (a)} If $G$ is connected, then $\Gamma_d(G) = \lceil
n/2 \rceil$. \\
\indent {\rm (b)} $\Gamma_d(G) \le 2n/3$ with equality if and only
if $G$ consists of disjoint copies of $K_3$. \label{p:cycle}
\end{prop}
\textbf{Proof.} (a) Suppose that $G$ is a cycle $C_n$. If $n$ is
even, $G$ has a perfect matching, and so, by
Theorem~\ref{t:uppbd}(c), $\Gamma_d(G) \le n/2$. If $n$ is odd, then
$\alpha'(G) = (n-1)/2$. By Theorem~\ref{t:uppbd}(b), $\Gamma_d(G)
\le n - \alpha'(G) = n - (n-1)/2 = (n+1)/2$. In both cases,
$\Gamma_d(G) \le \lceil n/2 \rceil$. To show that $\Gamma_d(G) \ge
\lceil n/2 \rceil$, we note that if $D$ is a directed cycle $C_n$,
then every vertex out-dominates itself and exactly one other vertex,
and so $\Gamma_d(G) \ge \gamma(D) = \lceil n/2 \rceil$. This proves
part~(a).

(b) To prove part~(b), let $G_1, G_2, \ldots, G_k$ be the components
of $G$, where $k \ge 1$. For $i = 1,2,\ldots,k$, let $G_i$ have
order~$n_i$. Since each component of a cycle, $n \ge 3k$. Applying
the result of part~(a) to each component of $G$, we have
\[
\Gamma_d(G) = \sum_{i=1}^k \Gamma_d(G_i) \le \sum_{i=1}^k \left(
\frac{n_i+1}{2} \right) = \frac{n + k}{2} \le  \frac{2n}{3},
\]
with equality if and only if $n = 3k$, i.e., if and only if $G_i =
C_3$ for each $i = 1,2,\ldots,k$.~$\Box$


\medskip
We remark that the upper bound of Theorem~\ref{t:regular1} can be
improved using tight lower bounds on the size of a maximum matching
in a regular graph established in~\cite{HeYe07}. Applying
Theorem~\ref{t:uppbd}(a) to these matching results in~\cite{HeYe07},
we have the following result. We remark that the $(n+1)/2$ bound in
the statement of Theorem~\ref{regular} is only included as it is
necessary when $n$ is very small or $r=2$.

\begin{thm}
For $r \ge 2$, if $G$ is a connected $r$-regular graph of order~$n$,
then
\[
\Gamma_d(G) \le \left\{ \begin{array}{ll} \displaystyle{
 \max \left\{  \left(
\frac{r^2 + 2r}{r^2 + r + 2} \right) \times \frac{n}{2},
\frac{n+1}{2} \right\} } &
\mbox{if $r$ is even} \\
& \\
\displaystyle{ \frac{(r^3+r^2-6r+2) \, n + 2r - 2}{2(r^3-3r)} } &
\mbox{if $r$ is odd} \\
\end{array} \right.
\]
\label{regular}
\end{thm}



\medskip
We close this section with the following observation. Graphs $G$
satisfying $\chi'(G) = \Delta(G)$ are called \emph{class~1} and
those with $\chi'(G) = \Delta(G) + 1$ are \emph{class~2}.

\begin{ob}
Let $G$ be an $r$-regular graph of order~$n$. Then the following holds. \\
{\rm (a)} If $G$ is of class~1, then $\Gamma_d(G) \le n/2$. \\
{\rm (b)} If $r \ge n/2$, then $\Gamma_d(G) \le \lceil n/2 \rceil$.
 \label{t:regular2}
\end{ob}
\textbf{Proof.} (a) Consider a $r$-edge coloring of $G$. The edges
in each color class form a perfect matching in $G$, and so, by
Theorem~\ref{t:uppbd}(c), $\Gamma_d(G) \le n/2$.

(b) If $n = 2$, then the result is immediate. Hence we may assume
that $n \ge 3$. By Dirac's theorem, $G$ is hamiltonian, and so
$\alpha'(G) \ge \lfloor n/2 \rfloor$. By Theorem~\ref{t:uppbd}(b),
$\Gamma_d(G) \le n - \alpha'(G) \le n - \lfloor n/2 \rfloor = \lceil
n/2 \rceil$.~$\Box$

\subsection{Outerplanar Graphs}

Let $\cOP_n$ denote the family of all maximal outerplanar graphs of
order~$n$. We define $\Mop(n) = \max \{ \Gamma_d(G) \}$
where the maximum is taken over all graphs $G \in \cOP_n$.

\begin{thm}
$\Mop(n) = \lceil n/2 \rceil$. \label{t:outer}
\end{thm}
\textbf{Proof.} Let $G \in \cOP_n$. Since every maximal outerplanar
graph is hamiltonian, we observe by Observation~\ref{o:spanning} and
Proposition~\ref{p:cycle}(a), that $\Gamma_d(G) \le \Gamma_d(C_n) =
\lceil n/2 \rceil$. Since this is true for an arbitrary graph $G$ in
$\cOP_n$, we have $\Mop(n) \le \lceil n/2 \rceil$. Hence it suffices
for us to prove that $\Mop(n) \ge \lceil n/2 \rceil$. If $n = 3$,
then by Observation~\ref{o:induced}, $\Gamma_d(G) \ge \Gamma_d(C_n)
= \lceil n/2 \rceil$, as desired. Hence we may assume that $n \ge
4$, for otherwise the desired result follows.

For $n \ge 4$ even, we take a directed cycle $\overrightarrow{C_n}$
on $n \ge 4$ vertices and a selected vertex $v$ on the cycle, and we
add arcs from every vertex $u$, where $u$ is neither the in-neighbor
nor the out-neighbor of $v$ on $\overrightarrow{C_n}$, to the vertex
$v$. The resulting orientation $D$ of the underlying maximal
outerplanar graph has $\gamma_d(D) = n/2$. Hence for $n \ge 4$ even,
we have $\Mop(n) = n/2$.

It remains for us to show that for $n \ge 5$ odd, $\Mop(n) =
(n+1)/2$. For $n \ge 5$ odd, we take a directed cycle
$\overrightarrow{C_n} \colon v_1v_2 \ldots v_n v_1$ on $n$ vertices.
We now add the arcs from $v_i$ to $v_1$ for all odd $i$, where $3
\le i \le n-2$, and we add the arcs from $v_1$ to $v_i$ for all even
$i$, where $4 \le i \le n-1$. Let $G$ denote the resulting
underlying maximal outerplanar graph and let $D$ denote the
resulting orientation of $D$. We now consider an arbitrary DDS $S$
in $D$.

Suppose first that $v_1 \in S$. In order to dominate the $(n-1)/2$
vertices $v_{2i+1}$, where $1 \le i \le (n-1)/2$, in $D$ we must
have that $|S \cap \{v_{2i},v_{2i+1}\}| \ge 1$ for all $i =
1,2,\ldots,(n-1)/2$. Hence in this case when $v_1 \in S$, we have
$|S| \ge (n+1)/2$.

Suppose next that $v_1 \notin S$. Then, $v_2 \in S$. In order to
dominate the $(n-3)/2$ vertices $v_{2i}$, where $2 \le i \le
(n-1)/2$, in $D$ we must have that $|S \cap \{v_{2i},v_{2i-1}\}| \ge
1$ for all $i = 2,\ldots,(n-1)/2$. In order to dominate $v_1$, there
is a vertex $v_j \in S$ for some odd $j$, where $3 \le j \le n$. Let
$j$ be the largest such odd subscript for which $v_j \in S$. If $j =
n$, then $v_n \in S$ and $|S| \ge (n+1)/2$, as desired. Hence we may
assume that $j < n$. In order to dominate the vertex $v_i$ for $i$
odd with $j < i \le n$, we must have $v_{i-1} \in S$. In particular,
we have that $v_{j+1} \in S$ to dominate $v_{j+2}$, implying that
$|S \cap \{v_{j},v_{j+1}\}| = 2$ while for $i$ odd where $i \ne j$
and $3 \le i \le n-2$, we have $|S \cap \{v_{i},v_{i+1}\}| \ge 1$,
implying that $|S| \ge (n+1)/2$.

In both cases, $|S| \ge (n+1)/2$. Since $S$ is an arbitrary DDS in
$D$, we have $\gamma(D) \ge (n+1)/2$. Hence, $\Gamma_d(G) \ge
(n+1)/2$, implying that $\Mop(n) = (n+1)/2$.~$\Box$

\subsection{Perfect Graphs}

Recall that a \emph{perfect graph} is a graph in which the chromatic
number of every induced subgraph equals the size of the largest
clique of that subgraph. Characterization of perfect graphs was a
longstanding open problem. The first breakthrough was due to Lovász
in 1972 who proved the Perfect Graph Theorem.

\begin{unnumbered}{Perfect Graph Theorem}
A graph is perfect if and only if its complement is perfect.
\end{unnumbered}

Let $\alpha \ge 1$ be an integer and let $\cG_\alpha$ be the class of
all graphs $G$ with $\alpha \ge \alpha(G)$. We are now in a position
to present an upper bound on the directed domination number of a
perfect graph in terms of its independence number.

\begin{thm}
If $G \in \cG_\alpha$ is a perfect graph of order $n \ge \alpha$,
then
\[
\Gamma_d(G) \le \alpha \log \left(  \lceil n/\alpha \rceil + 1 \right).
\]
\label{main_thm}
\end{thm}
\textbf{Proof.} By the Perfect Graph Theorem, the complement $\barG$
of $G$ is perfect. Hence, $\chi(\barG) = \omega(\barG) = \alpha(G)$.
The desired result now follows from Theorem~\ref{color}.~$\Box$

\section{Interplay between Transversals and Directed Domination}

In this section, we present upper bounds on the directed domination
number of a graph by demonstrating an interplay between the directed
domination number of a graph and the transversal number of a
hypergraph.
%
%
We shall need the following upper bounds on the transversal number
of a uniform hypergraph established by Alon~\cite{Al90} and
Chv\'{a}tal and McDiarmid~\cite{ChMc}. Applying probabilistic
arguments, Alon~\cite{Al90} showed the following result.

\begin{thm}{\rm (Alon~\cite{Al90})}
For\/ $k \ge 2$, if $H$ is a $k$-uniform hypergraph with\/ $n$
vertices and\/ $m$ edges, then $\tau(H) \le (m+n)(\ln k)/ k$.
\label{transAlon}
\end{thm}


\begin{thm}{\rm (Chv\'{a}tal, McDiarmid~\cite{ChMc})}
For\/ $k \ge 2$, if $H$ is a $k$-uniform hypergraphs with\/ $n$
vertices and\/ $m$ edges, then \/ $\TR{H} \le (n + \lfloor
\frac{k}{2} \rfloor m) / \lfloor \frac{3k}{2} \rfloor$.
bound is sharp. \label{ChMcbound}
\end{thm}

We proceed further with two lemmas. For this purpose, we shall need
the Szekeres-Wilf Theorem.

\begin{thm}{\rm (Szekeres-Wilf~\cite{SzWi68})}
If $G$ is a $k$-degenerate graph, then $\chi(G) \le k+1$.
\label{Wilf}
\end{thm}

\begin{lem}
If $G$ is a graph and $D$ is an orientation of $G$ such that
$\Delta^-(D) \le k$ for some fixed integer $k \ge 0$, then $\chi(G)
\le 2k+1$. \label{l:degen}
\end{lem}
\textbf{Proof.} It suffices to show that $G$ is $2k$-degenerate,
since then the desired result follows from the Szekeres-Wilf
Theorem. Assume, to the contrary, that $G$ is not $2k$-degenerate.
Then there is a subset $S$ of $V(G)$ such that the subgraph $G_S =
G[S]$ induced by $S$ has minimum degree at least~$2k+1$ and hence
contains at least~$(2k+1)|S|/2$ edges. Let $D_S = D[S]$ be the
orientation of $D$ induced by $S$. Since $\Delta^-(D) \le k$, we
have that $\Delta^-(D_S) \le k$ and
\[
k|S| \ge \sum_{v \in V(D_S)} d^-(v) = |E(G_S)| \ge (2k+1)|S|/2
>k|S|,
\]
a contradiction.~$\Box$

\begin{lem}
Let $D$ be an orientation of a graph $G$. If $G$ contains $n_k$
vertices with in-degree at most~$k$ in $D$ for some fixed integer $k
\ge 0$, then $n_k \le (2k+1)\alpha(G)$. \label{l:orient1}
\end{lem}
\textbf{Proof.} Let $V_k$ denote the set of all vertices of $G$ with
in-degree at most~$k$ in $D$, and so $n_k = |V_k|$. Let $G_k =
G[V_k]$ and let $D_k = D[V_k]$. Then, $D_k$ is an orientation of
$G_k$ such that $\Delta^-(D_k) \le k$, and so by
Lemma~\ref{l:degen}, $\chi(G_k) \le 2k+1$. Since every color class
of $G_k$ is an independent set, and therefore has cardinality at
most~$\alpha(G)$, we have that $n_k = |V_k| \le \chi(G_k)\alpha(G)
\le (2k+1)\alpha(G)$.~$\Box$

\medskip
Let $f(n,k)$, $g(n,k)$, and $h(n,k)$ be the functions of $n$ and $k$
defined as follows.
\[
\begin{array}{lcl} \1
f(n,k) & = &  \displaystyle{  2n\ln(k+2)/(k + 2) + (2k+1)\alpha(G) }  \\ %
\1
g(n,k) & = & \displaystyle{ n(k+2)/3k + 2(2k+1)\alpha(G)/3 } \\ %
\1 %
h(n,k) & = & \displaystyle{ n(k+1)/(3k-1)
 + 2k(2k+1)\alpha(G)/(3k-1) }
\end{array}
\]


\begin{thm}
If $G$ is a graph on $n$ vertices, then
\[
\Gamma_d(G) \le \left\{ \begin{array}{ll}
 \displaystyle{ \min_{k \ge 0} \{ f(n,k), g(n,k) \} } &
\mbox{if $k$ is even} \\
& \\ \displaystyle{ \min_{k \ge 1} \{ f(n,k), h(n,k) \} } & \mbox{if
$k$ is
odd} \\
\end{array} \right.
\]
\label{m:uppbd1}
\end{thm}
\textbf{Proof.} Let $D$ be an arbitrary orientation of the graph $G$
and let $k \ge 0$ be an arbitrary integer. Let $V_k$ denote the set
of all vertices of $G$ with in-degree at most~$k$ in $D$ and let
$n_k = |V_k|$. Let $V_{> k} = V(G) \setminus V_k$, and so all
vertices in $V_{> k}$ have in-degree at least~$k+1$ in $D$. Let
$H_{> k}$ be the hypergraph obtained from the CINH $H_D$ of $D$ by
deleting the $n_k$ edges corresponding to closed in-neighborhoods of
vertices in $V_k$. Each edge in $H_{> k}$ has size at least $k + 2$.

We now define the hypergraph $H$ as follows. For each edge $e_v$ in
$H_{> k}$ corresponding to the closed in-neighborhood of a vertex
$v$ in $V_{> k}$, let $e_v'$ consist of $v$ and exactly $k+1$
vertices from $N^{-}(v)$. Thus, $e_v' \subseteq e_v$ and $e_v'$ has
size~$k+2$. Let $H$ be the hypergraph obtained from $H_{> k}$ by
shrinking all edges $e_v$ of $H_{> k}$ to the edges $e_v'$. Then,
$H$ is a $(k+2)$-uniform hypergraph with $n$ vertices and $n - n_k$
edges.

Every transversal $T$ in $H$ contains a vertex from the closed
in-neighborhood of each vertex from the set $V_{> k}$ in $D$, and
therefore $T \cup V_k$ is a DDS in $D$. In particular, taking $T$ to
be a minimum transversal in $H$, we have that $\gamma(D) \le \tau(H)
+ n_k$. By Lemma~\ref{l:orient1}, $n_k \le (2k+1)\alpha(G)$.
Applying Theorem~\ref{transAlon} to the hypergraph $H$, we have that
\[
\tau(H) \le (n + n - n_k)\ln (k+2)/ (k+2) \le 2n\ln(k+2)/(k+2),
\]

\noindent and so $\gamma(D) \le \tau(H) + n_k \le 2n\ln(k+2)/(k+2) +
\alpha(G)(2k+1) = f(n,k)$. Applying Theorem~\ref{ChMcbound} to the
hypergraph $H$ for $k$ even, we have that
\[
\tau(H) \le (2n + k(n-n_k))/3k = n(k+2)/3k - n_k/3,
\]

\noindent and so $\gamma(D) \le \tau(H) + n_k \le n(k+2)/3k + 2n_k/3
\le n(k+2)/3k + 2(2k+1)\alpha(G)/3 = g(n,k)$. Thus for $k$ even, we
have that $\Gamma_d(G) \le \min \{ f(n,k), g(n,k) \}$.
Applying Theorem~\ref{ChMcbound} to the hypergraph $H$ for $k$ odd,
we have that
\[
\tau(H) \le (2n + (k-1)(n-n_k))/(3k -1) = n(k+1)/(3k-1) -
(k-1)n_k/(3k -1),
\]
and so $\gamma(D) \le \tau(H) + n_k \le n(k+1)/(3k-1) + 2k
n_k/(3k-1) \le n(k+1)/(3k-1) + 2k(2k+1)\alpha(G)/(3k-1) = h(n,k)$.
Thus for $k$ odd, we have that $\Gamma_d(G) \le \min \{ f(n,k),
h(n,k) \}$.~$\Box$

\medskip
Let $f_n(\alpha)$, $g_n(\alpha)$, and $h_n(\alpha)$ be the functions
of $n$ and $\alpha$ defined as follows.
\[
\begin{array}{lcl} \1
f_n(\alpha) & \doteq &  \displaystyle{  \sqrt{2n \alpha} \left( \ln
(\sqrt{2n/\alpha}\,) + 2 \right) - 2\alpha }  \\ %
\1 g_n(\alpha) & \doteq & \displaystyle{ \frac{1}{3}\left(n +
2\alpha + 4
\sqrt{ 2 n \alpha } \right) } \\ %
\1 %
h_n(\alpha) & \doteq & \displaystyle{ \frac{1}{3}\left(n +
\frac{14}{3}\alpha + \frac{ \sqrt{2\alpha} \, ( 27n + 20 \alpha)
}{3\sqrt{5\alpha + 6n}} \right) }
\end{array}
\]

As a consequence of Theorem~\ref{m:uppbd1}, we have the following
upper bound on the directed domination of a graph.

\begin{thm}
If $G$ is a graph on $n$ vertices with independence number $\alpha$,
then
\[
\displaystyle{ \Gamma_d(G) \le  \min \left\{ f_n(\alpha),
g_n(\alpha), h_n(\alpha) \right\} .}
\]
\label{thm:uppbd2}
\end{thm}
\textbf{Proof.} By Theorem~\ref{m:uppbd1}, we need to optimize the
functions $f(n,k), g(n,k)$ and $h(n,k)$ over $k$ to obtain an upper
bound on $\Gamma_d(G)$. To simplify the notation, let $\alpha =
\alpha(G)$. Optimizing the function $g(n,k)$ over $k$ (treating $n$
as fixed), we get $g(n,k) \le g_n(\alpha)$, while optimizing the
function $h(n,k)$ over $k$ (treating $n$ as fixed), we get $h(n,k)
\le h_n(\alpha)$. Optimization of the function $f(n,k)$ is
complicated. Hence to simplify the computations, we choose a value
$k^*$ for $k$ and show that $f(n,k^*) \le f_n(\alpha)$.
Suppose $\alpha \ge  n/2$. Then, $\alpha = cn$ with $1 \ge c \ge
1/2$. Substituting this into $f_n(\alpha)$ we get $f_n(\alpha) =
n\sqrt{2c}(\ln (2/c) +2) - 2cn = n \left( \sqrt{2c} ( \ln (2/c) +2)
- 2c \right) \ge n$, and so the inequality $\Gamma_d(G) \le
f_n(\alpha)$ holds trivially. Hence we may assume that $\alpha \le
n/2$. We now take $k = \sqrt{2n/\alpha} -2 \ge 0$. Substituting into
$f(n,k) = 2n\ln (k+2)/(k+2) + (2k+1)\alpha$, we get
\[
\begin{array}{lcl} \1
f(n,k) & = & 2n\ln (\sqrt{2n/\alpha}\,)/\sqrt{2n/\alpha} +
(2\sqrt{2n/\alpha} - 3)\alpha \\
\1 & = & \sqrt{2n \alpha} \ln (\sqrt{2n/\alpha}\,) +
2\alpha\sqrt{2n/\alpha}  - 3\alpha \\
\1 & = &  \sqrt{2n \alpha} \left( \ln (\sqrt{2n/\alpha}\,) + 2
\right) - 3\alpha \\
\1 & < &  f_n(\alpha),
\end{array}
\]
as desired.~$\Box$

\medskip
If every edge of a hypergraph $H$ has size at least~$r$, we define
an \emph{$r$-transversal} of $H$ to be a transversal $T$ such that
$|T \cap e| \ge r$ for every edge $e$ in $H$. The
\emph{$r$-transversal number} $\TRr{H}$ of $H$ is the minimum size
of an $r$-transversal in $H$. In particular, we note that $\TRo{H} =
\TR{H}$. For integers $k \ge r$ where $k \ge 2$ and $r \ge 1$, we
first establish general upper bounds on the $r$-transversal number
of a $k$-uniform hypergraph. Our next result generalizes that of
Theorem~\ref{transAlon} due to Alon~\cite{Al90}, as well as
generalizes results due to Caro~\cite{Ca90c}.

\begin{thm}
For integers $k \ge r$ where $k \ge 2$ and $r \ge 1$, let $H$ be a
$k$-uniform hypergraph with\/ $n$ vertices and\/ $m$ edges. Then,
$\TRr{H} \le n\ln k/k  + rm(2\ln k)^r/k$.
 \label{t:rtrans}
\end{thm}
\textbf{Proof.} Pick every vertex of $V(H)$ randomly with
probability~$p$ to be determined later but such that $(1-p) > 1/2$.
Let $X$ be the set of randomly picked vertices and let
$E_X$ be the set of edges of $E(H)$ whose intersection with $X$ is
at most~$r-1$. For every fixed edge~$e \in E(H)$, the probability
that $e$ is in $E_X$ is exactly
\begin{equation}
\Pr(e \in E_X) = \sum_{i=0}^{r-1} {k \choose i}p^{i}(1-p)^{k-i} =
(1-p)^{k} \sum_{i=0}^{r-1} {k \choose i} \left( \frac{p}{1-p}
\right)^{i}. \label{Eq1}
\end{equation}

We now choose $p = \ln k/k$. With this choice of $p$, we have that
$(1-p) > 1/2$. Hence, $1/(1-p)^i < 2^i$ for all $i \ge 1$. Since $1
- x \le e^{-x}$ for all $x \in R$, we note that $(1-p)^k \le e^{-pk}
= e^{-\ln k} = 1/k$. Substituting $p = \ln k/k$ into
Equation~(\ref{Eq1}) we therefore get

\[
\Pr(e \in E_X) \le \frac{1}{k}  \sum_{i=0}^{r-1} \frac{k^i}{i!}
\cdot \frac{p^i}{(1-p)^i} \le  \frac{1}{k} \sum_{i=0}^{r-1}
\frac{(2kp)^i}{i!} \le  \frac{1}{k} \sum_{i=0}^{r-1} (2\ln k)^i \le
\frac{1}{k} (2\ln k)^r,
\]

\noindent since $1+ q+q^2+ \cdots + q^{r-1} = (q^r - 1)/(q-1) \le
q^r$ for $q > 1$ and $r \ge 1$. For each edge $e \in E_X$, we add $r
- |e \cap X|$ (which is at most~$r$) vertices from $e \setminus X$
to a set $Y$. Then, $T = X \cup Y$ is a $r$-transversal in $H$ and
$|Y| \le r|E_X|$. By the linearity of expectation, $E(T) = E(X) +
E(Y) \le  E(X) + r E(E_X) = n\ln k/k  + rm(2\ln k)^r/k$.~$\Box$

\medskip
Using $r$-transversals in hypergraphs, we obtain the following bound
on the directed $r$-domination number of a graph.

\begin{thm}
For $r \ge 1$ an integer, if $G$ is a graph on $n$ vertices, then
\[
\Gamma_d(G,r) \le \min_{k \ge r} \, \left\{ (2k-1)\alpha(G) +  n\ln
(k+1)/(k+1) + rn(2\ln (k+1))^r/(k+1) \right\}.
\]
\label{t:rdom}
\end{thm}
\textbf{Proof.} Let $D$ be an arbitrary orientation of the graph $G$
and let $k \ge r$ be an arbitrary integer. Let $V_{<k}$ denote the
set of all vertices of $G$ with in-degree at most~$k-1$ in $D$ and
let $n_{<k} = |V_{<k}|$. Let $G_{<k}$ be the subgraph of $G$ induced
by the set $V_{<k}$ and let $D_{<k}$ be the orientation of $G_{<k}$
induced by $D$. Then, $\Delta^-(D_{<k}) \le k-1$, and so, by
Lemma~\ref{l:degen}, $\chi(G_{<k}) \le 2k-1$, implying that $n_{<k}
\le (2k-1)\alpha(G)$.

Let $V_{k} = V(G) \setminus V_{<k}$, and so all vertices in $V_{k}$
have in-degree at least~$k$ in $D$. Let $H_{k}$ be the hypergraph
obtained from the CINH $H_D$ of $D$ by deleting the $n_{<k}$ edges
corresponding to closed in-neighborhoods of vertices in $V_{<k}$.
Each edge in $H_{k}$ has size at least~$k + 1$.
We now define the hypergraph $H$ as follows. For each edge $e_v$ in
$H_{k}$ corresponding to the closed in-neighborhood of a vertex $v$
in $V_{k}$, let $e_v'$ consist of $v$ and exactly $k$ vertices from
$N^{-}(v)$. Thus, $e_v' \subseteq e_v$ and $e_v'$ has size~$k+1$.
Let $H$ be the hypergraph obtained from $H_{k}$ by shrinking all
edges $e_v$ of $H_{k}$ to the edges $e_v'$. Then, $H$ is a
$(k+1)$-uniform hypergraph with $n$ vertices and $n - n_{<k}$ edges.

Every $r$-transversal $T$ in $H$ contains at least~$r$ vertices from
the closed in-neighborhood of each vertex from the set $V_{k}$ in
$D$, and therefore $T \cup V_{<k}$ is a DrDS in $D$. In particular,
taking $T$ to be a minimum $r$-transversal in $H$, we have that
$\gamma_r(D) \le \tau_r(H) + n_{<k}$. By Lemma~\ref{l:orient1},
$n_{<k} \le (2k-1)\alpha(G)$. Noting that $k+1 \ge r+1 \ge 2$, we
can apply Theorem~\ref{t:rtrans} to the hypergraph $H$ yielding
$\tau_r(H) \le n\ln (k+1)/ (k+1) + r (n - n_{<k})(2\ln
(k+1))^r/(k+1)$,
%
and so $\gamma_r(D) \le \tau_r(H) + n_{<k} \le (2k-1)\alpha(G) +
n\ln (k+1)/ (k+1) + r n (2\ln (k+1))^r/(k+1)$.
Since this is true for every integer $k \ge r$, the desired upper
bound on $\Gamma_d(G,r)$ follows.~$\Box$

\section{Open Questions}

We close with a list of open questions and conjectures that we have
yet to settle. Let $\cR_n$ denote the family of all $r$-regular
graphs of order~$n$. We define $m(n,r) = \min \{ \Gamma_d(G)  \}$
and $M(n,r) = \max \{ \Gamma_d(G)  \}$, where the minimum and
maximum are taken over all graphs $G \in \cR_n$. Then, $m(n,1) =
M(n,1) = n/2$. By Proposition~\ref{p:cycle}, $m(n,2) = n/2$ while
$M(n,2) = 2n/3$. We remark that by Theorem~\ref{t:regular1}, for $r
\ge 2$, we know that
\begin{equation}
\frac{n}{2} \le M(n,r) \le \left( \frac{r+2}{r+1} \right) \cdot
\frac{n}{2} \label{EqM}
\end{equation}
(and this upper bound on $M(n,r)$ can be improved slightly by
Theorem~\ref{regular}).

\noindent \textbf{Conjecture~1.} For $r \ge 3$, $M(n,r) = n/2$.

By Theorem~\ref{t:lowerbd}(a), we know that if $G \in \cR_n$, then
$\Gamma_d(G) \ge \alpha(G) \ge n/(r+1)$, and so $n/(r+1) \le
m(n,r)$. Moreover taking $n/(r+1)$ copies of $K_{r+1}$, we have by
Theorem~\ref{t:KnLow} that $m(n,r) \le n \log (r+2)/(r+1)$. We pose
the following question.

\noindent \textbf{Question~1.} For $r \ge 3$, does there exists a
constant $c$ such that $m(n,r) \le cn/(r+1)$?

Let $\cOP_n$ denote the family of all maximal outerplanar graphs of
order~$n$ and define $\mop(n) = \min \{ \Gamma_d(G)  \}$, where the
minimum is taken over all graphs $G \in \cOP_n$. Since outerplanar
graphs are $3$-colorable, we note by Theorem~\ref{t:lowerbd}(b) that
for every graph $G \in \cOP_n$, $\Gamma_d(G) \ge n/3$, implying that
$\mop(n) \ge n/3$. By Theorem~\ref{t:outer}, we know that $\mop(n)
\le \lceil n/2 \rceil$. Thus, $n/3 \le \mop(n) \le \lceil n/2
\rceil$.

\noindent \textbf{Problem~1.} Find good lower and upper bounds on
$\mop(n)$.

Let $\cPP_n$ denote the family of all maximum planar graphs of
order~$n$. We define $\mpp(n) = \min \{ \Gamma_d(G)  \}$ and
$\Mpp(n) = \max \{ \Gamma_d(G) \}$, where the minimum and maximum
are taken over all graphs $G \in \cPP_n$.

\noindent \textbf{Problem~2.} Find good lower and upper bounds on
$\mpp(n)$ and $\Mpp(n)$.

\newpage

\end{document}